\documentclass[a4paper,12pt]{article} 
\usepackage{amsmath, amsthm, amssymb}
\usepackage{url}

\usepackage{physics}
\usepackage{graphicx,lipsum}
\graphicspath{ {./downloads/} }

\usepackage[colorlinks,citecolor=red,urlcolor=blue,bookmarks=false,hypertexnames=true]{hyperref}
\usepackage{tikz}
\usetikzlibrary{calc}
\usetikzlibrary{shapes}
\usepackage[autostyle]{csquotes}
\makeatletter

\usepackage[colorlinks]{hyperref}
\usepackage[nameinlink,capitalize]{cleveref}
\newtheorem{theorem}{Theorem}
\newtheorem{corollary}{Corollary}[theorem]

 \newtheorem{proposition}{Proposition}[section]
\newtheoremstyle{named}{}{}{\itshape}{}{\bfseries}{.}{.5em}{\thmnote{#3's }#1}
\theoremstyle{named}

\newtheorem{remark}{Remark}
\theoremstyle{definition}
\newtheorem{definition}{Definition}
\newtheorem{example}{Example}

\theoremstyle{remark}

\usepackage{xspace}
\usepackage[margin=1.0in]{geometry}

\usepackage{titlesec}

\usepackage{mathtools}

\DeclarePairedDelimiterX{\inp}[2]{\langle}{\rangle}{#1, #2}
\titleformat{\chapter}
  {\Large\bfseries} 
  {}                
  {0pt}            
  {\huge}

\begin{document}

\begin{center}
\fontsize{13pt}{10pt}\selectfont
    \textsc{\textbf{ PRIME PRINCIPAL RIGHT IDEAL RINGS}}
    \end{center}
\vspace{0.1cm}
\begin{center}
   \fontsize{12pt}{10pt}\selectfont
    \textsc{{\footnotesize $Tamem \ Al-Shorman^* and\  Malik\ Bataineh \ $ }}
\end{center}
\vspace{0.2cm}
\begin{center}
   \fontsize{12pt}{10pt}\selectfont
    \textsc{Abstract}
\end{center}
Let R be a commutative ring with unity $1\in R$. In this article, we introduce the concept of prime principal right ideal rings (\textbf{PPRIR}), A prime ideal P of R is said to be prime principal right ideal (\textbf{PPRI}) is given by $P =\{ ar : r\in R\}$ for some element a. The ring R is said to be prime principal right ideal ring (\textbf{PPRIR}) if every prime ideal of R is a prime principal right ideal (\textbf{PPRI}). A prime principal right ideal ring R is called a prime principal right ideal domain (\textbf{PPRID}) if R is a domain. Several properties and characteristics of prime principal right ideal ring (\textbf{PPRIR}). 

\begin{center}
   \fontsize{12pt}{10pt}\selectfont
    \textsc{{1. Introduction}}
\end{center}
Throughout this article, all rings are commutative ring with nonzero unity 1. An ideal P of the ring R is said to be prime ideal if for all $x,y \in R$ with $xy \in P$ then $x\in P$ or $y\in P$. Since prime ideals are so significant in commutative rings theory, this concept has been expanded and examined in a variety of ways. Some of these generalizations are as important as the prime ideals.  Burton in \cite{patterson1972first} defined the maximal ideals of R, an ideal P of R is said to maximal ideal provided that $P \neq R$ and whenever I is an ideal of R with $P \subset I \subseteq R$, then J = R. A prime ideal P is said to be a prime ideal minimal  over an ideal I if it is minimal among all prime ideals containing I. 
 
 Let P be a proper ideal of R. Then the radical of P is denoted by Rad(P) and it is defined as written below:
 \begin{center}
     Rad(P)=$\{ a\in R : a^n \in P $ for some positive integer n$\}$
 \end{center}
Note that Rad(P) is actually an ideal of the ring R which contains P (see \cite{patterson1972first}).

Let P be an ideal of R, P is said to be principal right ideal is given by $P = \{ar : r\in R\}$ for some element a. Conrad defined Noetherian rings in \cite{conradnoetherian}. A commutative ring R is called Noetherian if each ideal in R is finitely generated. In \cite{johnson1963principal}, R is said to be principal right ideal ring if every right ideal of R is a principal right ideal. We say that R is a principal right ideal domain if R is a domain.

A prime principal left ideal ring (\textbf{PPLIR}) and domain (\textbf{PPLID}) are defined symmetrically. A ring R is called a prime principal ideal ring (\textbf{PPIR}) if R is both a prime principal right ideal ring (\textbf{PPRIR}) and a prime principal left ideal ring (\textbf{PPLIR})

In Section Two, we introduce the concept of prime principal right ideal rings (\textbf{PPRIR}). We say that P is a prime principal right ideal (\textbf{PPRI}) of R is given by $P =\{ ar : r\in R\}$ for some element a, when P is a prime ideal of R. A ring R we say that a prime principal right ideal ring (\textbf{PPRIR}) if every prime ideal is prime principal right ideal. A prime principal right ideal ring R is called a prime principal right ideal domain (\textbf{PPRID}) if R is a domain. We study rings with this property. Clearly, every  principal right ideal rings (\textbf{PRIR}) is prime principal right ideal rings (\textbf{PPRIR}). In Theorem \ref{thm1}, we show R/P is a prime principal right ideal domain (\textbf{PPRID}) when R is a prime principal right ideal ring (\textbf{PPRIR}) and P is a prime ideal of R. We investigate some basic properties of prime principal right ideal ring (\textbf{PPRIR}).

\begin{center}
   \fontsize{12pt}{10pt}\selectfont
    \textsc{{2. Prime Principal Right Ideal Rings}}
\end{center}

In this section, we introduce and study prime principal right ideal ring (\textbf{PPRIR}).

\begin{definition}\label{def1}
Let R be a commutative ring with identity $1\in R$.
\\
1. A prime ideal P of R is prime principal right ideal (\textbf{PPRI}) is given by $P =\{ ar : r\in R\}$ for some element a.
\\
2. R is a prime principal right ideal ring (\textbf{PPRIR}) if every prime ideal is prime principal right ideal. A prime principal right ideal ring R is called a prime principal right ideal domain (\textbf{PPRID}) if R is a domain.  
 \end{definition}

 Clearly, every  principal right ideal rings (\textbf{PRIR}) is prime principal right ideal rings (\textbf{PPRIR}). We show examples of prime principal right ideal rings in the following examples. For a ring R and $n\in \mathbf{Z^+}$, $M_n(R)$ means the ring of $n\times n$ matrices over R and $I_n(R)=\{rI: r\in R\}$, where $I$ is the $n\times n$ identity matrix of $M_n(R)$.

 \begin{example}\label{ex1}
 1. Every field and division ring are prime principal right ideal ring.
 
 2. Let R be a principal right ideal ring and J a subring of $M_n(R)$ containing $I_n(R)$. Then J is a principal right ideal ring (see \cite{goodearl2004introduction},  Proposition 1.7), so J is a prime principal right ideal ring. In particular, $M_n(R)$ is a prime principal right ideal ring. 
 \end{example}
 
 \begin{theorem}\label{thm1}
 Let R be a prime principal right ideal ring (\textbf{PPRIR}). Then P is a prime ideal of R if and only if R/P is a prime principal right ideal domain (\textbf{PPRID}).
 \\
 \\
 \textbf{Proof.} First, take P is a prime ideal of R and R is (\textbf{PPRIR}), so is the quotient ring R/P. The only thing left is to make sure R/P doesn't have any zero divisors. Assume the following:
 \begin{center}
     (x + P)(y + P) = P
 \end{center}
 In other words,  the zero element of the ring R/P is the product of these two cosets. The preceding equation clearly equates to a requirement that xy + P = P, or $xy \in P$. But by assumption P is a prime ideal, so $x\in P$ or $y \in P$. This is means that either the coset  x + P = P or y + P = P. Hence, R/P is without zero divisors. Therefore R/P is a (\textbf{PPRID}).
 
 We just invert the argument to establish the converse. Assume that R/P is (\textbf{PPRID}) and $xy \in P$. This indicates, in terms of cosets, that
 \begin{center}
     (x + P)(y + P) = xy + P = P
 \end{center}
 By assumption R/P has not zero divisors, so x + P = P or y + P = P. In each case, one of x or y belongs to p, forcing P to be a prime ideal of R.
 \end{theorem}

Note that, let R be a commutative ring with unity then every maximal ideal is a prime ideal (see \cite{patterson1972first}).

\begin{proposition}\label{prop1}
If R is a prime principal right ideal ring (\textbf{PPRIR}), then every maximal ideal is prime principal right ideal (\textbf{PPRI}). 
\\
\\
\textbf{Proof.} Assume that R is a prime principal right ideal ring (\textbf{PPRIR}), so R is a commutative ring with unity. Then every maximal ideal is prime ideal, but every prime ideal of R is a prime principal right ideal (\textbf{PPRI}). Hence the maximal ideal is a prime principal right ideal (\textbf{PPRI}).
 
\end{proposition}
 
The converse of Proposition \ref{prop1} does not hold, the following example show that there exist prime ideal of the prime principal right ideal ring (\textbf{PPRIR}) which is not maximal ideal.

\begin{example}\label{ex2}
Let $R=\mathbf{Z} \times \mathbf{Z} $ be a prime principal right ideal ring and $P=\mathbf{Z} \times \{0\}$ be a prime ideal of R. Since
\begin{center}
    $\mathbf{Z} \times \{0\} \subset \mathbf{Z} \times \mathbf{Z}_e \subset R $
\end{center}
 with $\mathbf{Z} \times \mathbf{Z}_e$ an ideal of R, then $\mathbf{Z} \times \{0\}$ is not maximal ideal.
\end{example}

\begin{proposition}\label{prop2}
Let R be a Boolean ring and (\textbf{PPRIR}). Then P is a prime principal right ideal (\textbf{PPRI}) of R if and only if P is a maximal ideal of R.
\\
\\
\textbf{Proof.} It suffices to show that if P is the prime principal right ideal (\textbf{PPRI}), P is also the maximal ideal. Suppose that I is an ideal of R with $P \subset I \subseteq R$, what we need to show is that $I = R$. Let $x\in I$ and $x\notin P$, then $x(1 - x) = 0 \in P$. But P is a prime principal right ideal of R with $x \notin P$, so $1 - x \in P \subset I$. Then x and 1 - x in I, it follows that $1 = x + (1 - x) \in I$. As a result, the ideal I contains the identity and, as a result, $I = R$. Since there is no appropriate ideal between P and the entire ring R, we come to the conclusion that P is a maximal ideal.
\end{proposition}

\begin{proposition}\label{prop3}
Let R be a prime principal right ideal ring (\textbf{PPRIR}). If P is a prime principal right ideal (\textbf{PPRI}) of R, then P id a semiprime ideal of R.
\\
\\
\textbf{Proof.} Assume that P is a prime principal right  (\textbf{PPRI}) of R, so R/P is a prime principal right ideal domain (\textbf{PPRID}) by Theorem \ref{thm1}, and so R/P has no zero divisors. Thus R/P has no nonzero nilpotent elements. Hence P is a semiprime ideal of R.  
\end{proposition}

\begin{theorem}\label{thm2}
A ring R is a prime principal right ideal ring  (\textbf{PPRIR}) if and only if for each infinite increasing sequence of prime ideals $P_1 \subset P_2 \subset P_3 \subset ... $ in R, then $P_n = P_{n+1}$ for all large n.
\\
\\
\textbf{Proof.} ($\Rightarrow$) If $P_1 \subset P_2 \subset P_3 \subset ...$ is an increasing sequence of prime ideals. Let $P = \bigcup\limits_{n \geq 1}^{} P_n$ be a prime ideal of R. Since R is a prime principal right ideal ring (\textbf{PPRIR}), so P is a Prime principal right ideal (\textbf{PPRI}). Now, by increasing condition for each finite subset of P lies in a common $P_n$, so a finite set of P is in some $P_m$. Thus $P\subset P_m$ as well as that $P_m \subset P$, so $P = P_m$. Then for all $n \geq m$, $ P_m \subset P_n \subset P = P_m$, so $P_n = P_m$.
\\
\\
($\Leftarrow$) Suppose that R is a not prime principal right ideal ring, so R has a prime ideal P that is not prime principal right ideal. Take $p_1 \in P$. Since P is not prime principal right ideal, so $P \neq \{p_1r : r\in R \}= (p_1)$. Consequently, there is a $p_2 \in P - (p_1)$ but $P \neq (p_1 , p_2)$. Follow the same steps to take $p_n$ in P for all $n \geq 1$ by making $p_n \in P -(p_1 , p_2 , ... , p_{n-1})$ for all $n \geq 2$. Then we have an increasing sequence of prime ideals $(p_1) \subset (p_1 , p_2) \subset ... \subset (p_1 , p_2 , ... , p_{n-1}) \subset ...$ in R where each prime ideal is completely contained within the one after it, so this is contradiction. Hence R is a prime principal right ideal ring (\textbf{PPRIR}).
\end{theorem}

\begin{theorem}\label{thm3}
Let R be a prime principal right ideal domain (\textbf{PPRID}). If $\phi$ from R to R is a surjective ring homomorphism, then $\phi$ is an injective ring homomorphism and thus is an isomorphic ring homomorphism.
\\
\\
\textbf{Proof.} Let $\phi : R \rightarrow R$ be a surjective ring homomorphism. The nth iteration $\phi ^n$, let $P_n = ker(\phi ^n)$ be a prime ideal of R and these prime ideals form an increasing chain $P_1 \subset P_2 \subset P_3 \subset ...$ since  $x \in P_n$ then $\phi^n(x)=0$ and then $\phi^{n+1}(x)=\phi(\phi^n(x))=\phi(0)=0 $, so $x \in P_{n+1}$. Since R is a prime principal right ideal ring (\textbf{PPRIR}) by Theorem \ref{thm2} $P_n = P_{n+1}$ for some n. Take $x \in ker(\phi)$, so $\phi(x)=0$. The map $\phi^n$ is surjective map, so $x = \phi^n(y)$ for some $y\in R$. Thus $0= \phi(x) = \phi(\phi^n(y))= \phi^{n+1}(y)$.Therefore $x\in ker(\phi^{n+1})= ker(\phi^n)$, so $x= \phi^n(y)=0$. Thus $ker(\phi)=\{0\}$, so $\phi$ is an injective ring homomorphism. 
\end{theorem}

The converse of Theorem \ref{thm3} is not true, the following example show this.
\begin{example}\label{ex3}
Let R[X] be a prime principal right ideal domain (\textbf{PPRID}) since it is a principal right ideal domain (\textbf{PRID}) and let $\phi: R[X] \rightarrow R[X]$ by $\phi(x)=x^2$. The homomorphism ring $\phi$ is an injective but not surjective.
\end{example}

\begin{theorem}\label{thm4}
If R is a prime principal right ideal ring (\textbf{PPRIR}) then R[X] is a prime principal right ideal ring (\textbf{PPRIR}).
\\
\\
\textbf{Proof.} It is clearly if $R = \{0\}$, so suppose that $R \neq \{0\}$. To prove the theorem must prove every prime ideal of R[X] is a prime principal right ideal (\textbf{PPRI}). By contradiction, assume that P is a prime ideal of R[X] bu not prime principal right ideal. We have $P \neq (0)$. Defined a sequence of polynomials $f_1 , f_2, f_3, ...$ in P as the following:
\\
1. Take $f_1 \in P-(0)$ with minimal degree. 
\\
2.Since $P \neq (f_1)$, take $f_2 \in P-(f_1)$ with minimal degree. Note that $deg(f_1) \leq deg(f_2)$.
\\
3. For $n\geq 2$, if we have defined $f_1 , f_2, ... , f_n$ in P then $P \neq (f_1 , f_2, ... , f_n)$,  so we may take $f_{n+1} \in P - (f_1 , f_2, ... , f_n)$ with minimal degree. 

We have $deg(f_n) \leq deg(f_{n+1})$ for all n. The case $n = 1$ was already examined and for $n \geq 2$, $f_n$ and $f_{n+1}$ are in $P - f_1 , f_2, ... , f_{n-1}$ so $deg(f_n) \leq deg(f_{n+1})$.

Let $d_n = deg(f_n)$ for $n \geq 1$ and $c_n$ be the leading coefficient of $f_n$, so $d_n \leq d_{n+1}$ and $f_n(X)= C_n X^{d_n} +$ terms of a lower degree.

Now choose the prime ideal $(c_1, c_2 , ...)$ in R is a prime principal right ideal (\textbf{PPRI}) since R is a prime principal right ideal ring (\textbf{PPRIR}). This prime ideal has elements are all linear combination of finitely many $c_n$, so $(c_1,c_2,...)= (c_1,c_2,...,c_m)$ for some m.

Since $c_{m+1}\in (c_1,c_2,...,c_m) $, we have  \[ c_{m+1}=  \sum_{n=1}^{m} x_n c_n\] for some $x_n \in R$. From $d_n \leq d_{m+1}$ for $n\leq m$, the leading term in $f_n(X) = c_n X^{d_n}+...$ can be transformed into a degree $d_{m+1}$ by using $f_n(X)X^{d_{m+1}-d_n}= c_n X^{d_{m+1}}+...$, and this is in P since $f_n(X) \in P$ and P is a prime ideal of R[X]. By the previous linear combination of $c_{m+1}$ the linear combination \[ \sum_{n=1}^{m} x_n f_n(X) X^{d_{m+1}-d_n} \] is in the ideal $(f_1,f_2,...,f_m)$ and its coefficient of $X^{d_{m+1}}$ is $\sum_{n=1}^{m} x_nc_n$ which equals the leading coefficient $c_{m+1}$ of $f_{m+1}(X)$ in degree $d_{m+1}$. The difference \[ f_{m+1}(X) - \sum_{n=1}^{m} x_n f_n(X) X^{d_{m+1}-d_n} \] is in P and not 0 since $f_{m+1} \in P- (f_1 , f_2, ... , f_m)$ with degree less than $d_{m+1}$. But $f_{m+1}$ has a minimal degree in $P- (f_1 , f_2, ... , f_m)$ and the previous difference is in $P- (f_1 , f_2, ... , f_m)$ with lower degree than $d_{m+1}$. This is the contradiction. Thus P is a prime principal right ideal (\textbf{PPRI}) of R[X]. 
\end{theorem}

In a single sentence, we would want to summarize this proof, " use a prime ideal of leading coefficients".

In the proof of Theorem \ref{thm4}, the prime principal right ideal ring (\textbf{PPRIR}) property of R where we mentioned is used $(c_1,c_2, c_3, ...) = (c_1,c_2,...,c_m)$ for some m. To get the contradiction in the evidence, all we need to do is $c_{m+1} \in (c_1,c_2,...,c_m)$ for some m. Since $(c_1) \subset (c_1,c_2) \subset (c_1, c_2, c_3) \subset ... $, the following property is what we require: for each infinite increasing sequence of prime ideals $P_1 \subset P_2 \subset P_3 \subset ...$ in R, $P_m=P_{m+1}$ for some m. Of course this is implied by prime principal right ideal ring (\textbf{PPRIR}) property, but it also implies the prime principal right ideal ring (\textbf{PPRIR}) property since a non-prime principal right ideal ring has an infinite increasing sequence of ideals with strict containment at each step see the proof of Theorem \ref{thm2}.

Where did we apply the assumption that R is a prime principal right ideal ring (\textbf{PPRIR}) in the proof of the Theorem \ref{thm4}? It is how we know the prime ideals $(c_1,c_2,...,c_n)$ for $n \geq 1$ stabilize for large n, so $c_{m+1} \in (c_1,c_2,...,c_m)$ for some m. The contradiction we obtain from that really shows $c_{m+1} \notin (c_1,c_2,...,c_m)$ for all m, so the proof of Theorem \ref{thm4} could be construed as proving the converse:  if R[X] is not prime principal right ideal ring then R is not prime principal right ideal ring.

\begin{remark}\label{rem1}
The converse of Theorem \ref{thm4} is true, if the ring R[X] is a prime principal right ideal domain (\textbf{PPRID}) then R is a prime principal right ideal domain (\textbf{PPRID}), since $R \cong R[x]/ (x)$.  
\end{remark}

\begin{corollary}\label{coro1}
If R is a prime principal right ideal ring then $R[X_1,...,X_n]$ is a prime principal right ideal ring.
\\
\\
\textbf{Proof.} The case of n=1 is in Theorem \ref{thm4}. For $n\geq 2$, write $R[X_1,...,X_n]$ as $R[X_1,...,X_{n-1}][X_n]$, with $R[X_1,...,X_{n-1}]$ being a prime principal right ideal ring  by the inductive hypothesis. As a result, we have been limited to the most basic case.
\end{corollary}

\begin{proposition}\label{prop4}
Let R be a prime principal right ideal ring (\textbf{PPRIR}). If P is a primary ideal of R then Rad(P) is a prime principal right ideal (\textbf{PPRI}) of R.
\\
\\
\textbf{Proof.} Suppose that P is a primary ideal of R, by proposition 4.1 in \cite{atiyah2018introduction} the Rad(P) is a prime ideal of R. But R is a prime principal right ideal ring (\textbf{PPRIR}), so Rad(P) is a prime principal right ideal of R. 
\end{proposition}

\begin{proposition}\label{}
Let R be a prime principal right ideal ring (\textbf{PPRIR}). If P is a prime ideal of R then Rad(P) is a prime principal right ideal (\textbf{PPRI}) of R.
\\
\\
\textbf{Proof.} Suppose that P is a prime ideal of R, so Rad(P) = P. But R is a prime principal right ideal ring (\textbf{PPRIR}), then every prime ideal of R is a prime principal right ideal (\textbf{PPRI}), and then Rad(P) is prime principal right ideal (\textbf{PPRI}). 
\end{proposition}

\begin{theorem}\label{thm5}
A ring R is a prime principal right ideal ring (\textbf{PPRIR}) if and only if every nonempty collection T of a prime ideals of R contains the maximum element in terms of inclusion there is a prime ideal in T not strictly contained in another prime ideal in T.
\\
\\
\textbf{Proof.} ($\Rightarrow$) Suppose that there is a nonempty collection T of a prime ideals in R which there is no maximum member in terms of inclusion. Thus suppose we begin with a prime ideal $P_1$ in R. we can find prime ideals recursively $P_2, P_3,...$ such that $P_{n-1} \subset P_n$ for all $n \geq 2$. $P_{n-1}$ would be a maximum element of T if there was no prime ideal in T strictly containing $P_{n-1}$, which does not exist. This is contradiction of Theorem \ref{thm2} .
\\
($\Leftarrow$) Let P be a prime ideal in R and T be the set of all a prime principal right ideal (\textbf{PPRI}) contained in P. Suppose that $I\in T$ no other element of T contains this, so I is a prime principal right ideal (\textbf{PPRI}) in P and no other prime principal right ideal (\textbf{PPRI}) of R contains I. we will prove P = I by contradiction, which would prove P is a prime principal right ideal (\textbf{PPRI}). If $P \neq I$, take $x \in P-I$, since I is a prime principal right ideal (\textbf{PPRI}), also $I + (x)$ is a prime principal right ideal (\textbf{PPRI}) where (x) is a prime ideal in R. So $I+(x) \in T$ and so $I \in I+(x)$ this is a contradiction of maximality of I as a member of T. Hence P = I.
\end{theorem}

\begin{theorem}\label{thm6}
Let R be a prime principal right ideal ring (\textbf{PPRIR}) and I be a proper ideal of R. If every prime ideal minimal over I is prime principal right ideal (\textbf{PPRI}), then there are only finitely many prime ideals minimal over I.
\\
\\
\textbf{Proof.} Let $T= \{ P_1 , ... P_n : $ each $P_i$ is a prime ideal minimal over I$\}$.  If $H\in T$ we have $H \subseteq I$, then any prime ideal P minimal over I contains some $P_i$, so $\{ P_1, ... , P_n\}$  is the set of minimal prime ideals of I. Now assume that $H \not\subseteq P$ for some $H \in T$. Let the set $S = \{J : J$ is a prime ideal of R with $I \subseteq J$ and $H \not \subseteq J$ for some $H \in T \}$. Since every element in T is a prime principal right ideal (\textbf{PPRI}), S is inductive and hence by Zorn's Lemma has a maximal element Q. But $I \subseteq Q$, Q contains a prime ideal $P$ minimal over I (see \cite{kaplansky1974topics}, Theorem 10). Thus $P \in T$, this is a contradiction.

\end{theorem}

\begin{center}
   \fontsize{12pt}{10pt}\selectfont
    \textsc{{3. Conclusions}}
\end{center}
In this study, we introduced the concept of prime principal right ideal rings which is a generalization of principal right ideal rings.  We investigated some basic properties of prime principal right ideal rigs. As a proposal to further the work on the topic, we are going to study the concepts of  S-prime principal right ideal rings and graded S-prime principal right ideal rings.

\bibliographystyle{amsplain}

\end{document}